\documentclass[12pt]{article}
\usepackage{latexsym}
\newsavebox{\toy}
\savebox{\toy}{\framebox[0.65em]{\rule{0cm}{1ex}}}
\newcommand{\QED}{\usebox{\toy}\end{demo}}
\newenvironment{property}%
{\begin{list}{}{\setlength{\rightmargin}{0pt}%
\setlength{\itemsep}{0pt}}}{\end{list}}
\newlength{\templength}
\newcommand{\bp}{\setlength{\templength}{\labelwidth}%
\setlength{\labelwidth}{2em}\begin{property}}
\newcommand{\ep}{\end{property}\setlength{\labelwidth}{\templength}}
\newtheorem{theorem}{\indent Theorem}[section]
\newtheorem{lemma}[theorem]{\indent Lemma}
\newtheorem{proposition}[theorem]{\indent Proposition}
\newtheorem{corollary}[theorem]{\indent Corollary}
\newtheorem{assumption}{\indent Assumption}
\newtheorem{definition}{\indent Definition}[section]
\newtheorem{remark}{\indent Remark}[section]
\newtheorem{exercise}{\indent Exercise}[section]
\newtheorem{example}{\indent Example}[section]
\newcommand{\Thm}[1]{Theorem \ref{Thm.#1}}
\newcommand{\Lem}[1]{Lemma \ref{Lem.#1}}

\newcommand{\Prop}[1]{Proposition \ref{Prop.#1}}

\newcommand{\Theorem}[1]{\begin{theorem}\label{Thm.#1}}
\newcommand{\Lemma}[1]{\begin{lemma}\label{Lem.#1}}
\newcommand{\Proposition}[1]{\begin{proposition}\label{Prop.#1}}
\newcommand{\Corollary}[1]{\begin{corollary}\label{Cor.#1}}
\newcommand{\Assumption}[1]{\begin{assumption}\label{Ass.#1}\rm}
\newcommand{\Definition}[1]{\begin{definition}\label{Def.#1}\rm}
\newcommand{\Remark}[1]{\begin{remark}\label{Rem.#1}\rm }
\newcommand{\Exercise}[1]{\begin{exercise}\label{Exe.#1}\rm }
\newcommand{\Example}[1]{\begin{example}\label{Exa.#1}\rm }
\newcommand{\bd}{\begin{displaymath}}
\newcommand{\ed}{\end{displaymath}}
\newcommand{\bdn}{\begin{equation}}
\newcommand{\bdnl}{\begin{equation}\label}
\newcommand{\edn}{\end{equation}}
\newcommand{\barray}{\begin{array}}
\newcommand{\earray}{\end{array}}
\newcommand{\bds}{\begin{description}}
\newcommand{\eds}{\end{description}}
\newcommand{\bitemize}{\begin{itemize}}
\newcommand{\eitemize}{\end{itemize}}
\newcommand{\benumerate}{\begin{enumerate}}
\newcommand{\eenumerate}{\end{enumerate}}
\newcommand{\btabbing}{\begin{tabbing}}
\newcommand{\etabbing}{\end{tabbing}}
\newcommand{\bcenter}{\begin{center}}
\newcommand{\ecenter}{\end{center}}
\newcommand{\bflushright}{\begin{flushright}}
\newcommand{\bflushleft}{\begin{flushleft}}
\newcommand{\eflushright}{\end{flushright}}
\newcommand{\eflushleft}{\end{flushleft}}
\newcommand{\bdnn }{\begin{eqnarray*}}
\newcommand{\ednn }{\end{eqnarray*}}
\newcommand{\bdmn}{\begin{eqnarray}}
\newcommand{\edmn}{\end{eqnarray}}

\newcommand{\nn}{\nonumber}
\newcommand{\SSC}[1]{\section{#1}\setcounter{equation}{0}}

\newcounter{biblio}
\newenvironment{references}%
{\begin{list}{[\arabic{biblio}]}{\usecounter{biblio}%
\setlength{\leftmargin}{2.5em}\setlength{\rightmargin}{0pt}%
\setlength{\labelwidth}{2em}\setlength{\itemsep}{0pt}}}{\end{list}}
\newcommand{\References}%
{\vspace{2.8ex plus .3ex minus .3ex}%
\begin{center}{\bf References}\end{center}\begin{references}}

\usepackage{amssymb}
\newcommand{\N}{{\mathbb{N}}}
\newcommand{\Z}{{\mathbb{Z}}}
\newcommand{\zd}{\Z^d}

\newcommand{\R}{{\mathbb{R}}}
\newcommand{\rd}{\R^d}
\newcommand{\C}{{\mathbb{C}}}
\newcommand{\ra }{\rightarrow }
\newcommand{\lra }{\longrightarrow }

\newcommand{\Lra}{\Longrightarrow }

\newcommand{\Llra}{\Longleftrightarrow }

\newcommand{\ov}{\overline}
\newcommand{\tl}{\widetilde}

\newcommand{\limn}{\lim_{n \nearrow \8}}

\newcommand{\vvs}{\vspace{2ex}}
\newcommand{\vs}{\vspace{1ex}}

\newcommand{\lan}{\langle \:}
\newcommand{\ran}{\: \rangle}
\newcommand{\lef}{\left}
\newcommand{\rig}{\right}

\newcommand{\st}{\stackrel}

\newcommand{\8}{\infty}

\newcommand{\6}{\partial}
\newcommand{\dps}{\displaystyle}
\newcommand{\sub}{\subset}

\newcommand{\bsh}{\backslash}

\renewcommand{\Re}{\mathop{\rm Re}\nolimits}

\renewcommand{\b}{\beta}
\newcommand{\gm}{\gamma}

\newcommand{\del}{\delta}
\newcommand{\D}{\Delta}

\newcommand{\e}{\varepsilon}

\newcommand{\h}{\eta}
\newcommand{\tht}{\theta}

\newcommand{\lm}{\lambda}

\newcommand{\m}{\mu}

\newcommand{\rh}{\rho}

\newcommand{\s}{\sigma}

\renewcommand{\t}{\tau}

\newcommand{\vp}{\varphi}

\newcommand{\w}{\omega}
\newcommand{\om}{\omega}
\newcommand{\W}{\Omega}

\newcommand{\cC }{{\cal C}}

\newcommand{\cF }{{\cal F}}
\newcommand{\cG }{{\cal G}}

\newcommand{\cN }{{\cal N}}
\newcommand{\cO }{{\cal O}}

\makeatletter
\def\section{\@startsection{section}{1}{\z@}{-3.5ex plus -1ex minus 
 -.2ex}{2.3ex plus .2ex}{\bf}}
\def\subsection{\@startsection{subsection}{2}{\z@}{-3.25ex plus -1ex minus 
 -.2ex}{1.5ex plus .2ex}{\bf}}
\makeatother

\newcommand{\bz}{\bar \zeta}
\newcommand{\md}{\mu^{(2)}_\8}

\newcommand{\IW}{{\mathbb{W}}}

\begin{document}

\bcenter 

\large{\bf Directed Polymers in 
Random Environment are Diffusive 
at Weak Disorder}\footnote{Submited November 2004, revised May 2005}

\vvs

\vvs \normalsize

\noindent Francis COMETS\footnote{Corresponding author. 
 Partially supported by CNRS, UMR 7599 ``Probabilit\'es et Mod\`eles 
al\'eatoires''\\
{\tt http://www.proba.jussieu.fr/pageperso/comets/comets.html}
}
%%%\footnote{Corresponding author}
\\

\vs \small
Universit{\'e} Paris 7, \\
Math{\'e}matiques, Case 7012\\
 2 place Jussieu,
75251 Paris, France \\
{\tt email: comets@math.jussieu.fr} \\

\vvs \normalsize

\noindent Nobuo YOSHIDA\footnote{Partially supported 
by JSPS Grant-in-Aid for Scientific
Research, Wakatekenkyuu (B) 14740071
{\tt http://www.math.kyoto-u.ac.jp/}$\widetilde{}$ {\tt nobuo/}}\\

\vs \small 
Division of Mathematics \\
Graduate School of Science \\
Kyoto University,\\
Kyoto 606-8502, Japan.\\
{\tt email: nobuo@math.kyoto-u.ac.jp}\\

\ecenter
\normalsize 

\begin{abstract}
In this paper, we consider directed polymers in random environment
with discrete space and time. 
For transverse dimension at least equal to
3, we prove that diffusivity holds for the path
in the full weak disorder region, i.e., where the 
partition function differs from its annealed value only by a non-vanishing 
factor. Deep inside this region, we also show that the quenched averaged 
energy has fluctuations of order 1. In complete generality
(arbitrary dimension and temperature), we prove
monotonicity of the phase diagram in the temperature.

\end{abstract}
\vspace{1cm}
\footnotesize 
\noindent{\bf Short Title.} Diffusive Directed Polymers
 


\normalsize
\newpage

%%%%%%%%%%%%%%%%%%%%
%
%\footnotesize 
%
%\tableofcontents
%\vspace{1cm}
%\normalsize
%%%%%%%%%%%%
%%%%%%%%%%%%%%
\vs
\normalsize
%%%%%%%%%%%%
\SSC{Introduction} \label{Intro}
%%%%%%%%%%%%%%%%%
In this classical model, the polymer is a long chain of size $n$ 
in the $1+d$-dimensional space, which is directed: It
stretches in the first direction of $\Z^{1+d}$, and therefore is
modelled as a graph $\{(t,\w_t) \}_{t=1}^n$, where 
$\w =(\w_t )_{t \in \N}$ is a nearest neighbor path in $\Z^d$. 
We introduce the probability space $(\W, \cF, P)$, which consists 
of the set $\W$ of all nearest neighbor paths in $\zd$, 
the cylindrical $\s$-field $\cF$, and the distribution $P$ 
of the $d$-dimensional simple random walk with $\w_0=0$. 
On the other hand, the environment describes locations 
which can be favorable or hostile
to the monomers: it is given by  
independent identically distributed random variables 
$\h=\{\h(n,x); n \in \N, x \in \Z^d\}$ with all finite exponential moments, 
defined on a probability space $(H, \cG, Q)$.
The polymer is attracted
by large positive values of the environment, and repelled by large negative 
ones. Further motivations for the model can be found in the physics 
literature
\cite{FiHu91}, \cite{KrSp91}, and a rigorous survey in \cite{CSY04}.
All these ingredients are incorporated in the {\bf polymer measure}
with environment $\h$:
\bdnl{mu_n}
\mu_n (d\om)= Z_n^{-1} \exp\{ \b H_n(\om)\} \; P(d\om)\;,
\edn
with
$$
H_n(\om)=\sum_{t=1}^n \h(t, \om_t).
$$
%%%%%
%\;,\quad Z_n=P[\exp\{ \b H_n(\om)\}]\;.
%%%%%%%%
%Here, $P$ denotes the distribution of the simple random walk 
%on the integer lattice $\Z^d$, 
%%%%%%%%
Here, $\b >0$ denotes the ``temperature inverse''
and prescribes how strongly the polymer
path $\om$ interacts with the medium, and the ``partition function''
$Z_n=P[\exp\{ \b H_n(\om)\}]$ is the normalizing constant 
making $\mu_n$ a probability measure on the path space. 
Here, and in the sequel, $P[X]$ stands for the $P$-expectation 
of a random variable $X$ on $(\W, \cF, P)$. 
%%%%%%%
%The reader 
%will make the distinction between $P$ and the law $Q$ of the  environment 
%$\h$. 
%%%%%%%%%%%
Note that the measure
$\mu_n$ depends on $n, \b$ and on the environment $\h$.
We denote by $\lm$ the function
\bdnl{lm(beta)}
\lm(\b)=\ln Q[ \exp\{\b \h(t,x)\}] \in \R\;,\quad \b \in \R.
\edn
Consistently with the notation $P[X]$, $Q[Y]$ stands for the $Q$-expectation 
of a random variable $Y$ on $(H, \cG, Q)$. We  assume that 
$\lm (\cdot )$ is finite on the whole real line. 
%%%%%%%%%%%
\Remark{b<0}
%%%%%%%%%%
%%%%We can naturally define the polymer measure (\ref{mu_n}) with negative 
%%%%$\b$, not only with positive ones. However, considering negative $\b$ 
%merely amounts to considering $-\h$ with $|\b|$ as the inverse temperature. 
%For this reason, we can restrict ourselves to positive $\b$ without loss, 
%as far as real $\b$ is concerned.
%Moreover, this restriction helps us simplify the exposition of some results 
%%%%in this paper, e.g., \Thm{beta_c} below.  
The definition (\ref{mu_n}) makes perfect sense for real 
$\b$'s. However, considering negative $\b$ 
merely amounts to considering $-\h$ with $|\b|$ as the inverse temperature. 
Without loss of generality, 
we will restrict ourselves to positive $\b$, 
which will help us simplify the statements of some results 
in this paper, e.g., \Thm{beta_c} below. 
%%%%%%%%%%%
\end{remark}
%%%%%%%%%%%%%%
The issue is to understand the asymptotics of the polymer $\om$ as $n \to \8$
under the  measure $\mu_n$, for typical realization of the environment.
In particular, one would like to determine the exponent $\xi=\xi (d, \b)
 \in [1/2,1)$ such that
$$ |\om_n| {\rm \ is \ of\ order \ } n^\xi$$
as $n \to \8$. 
Another --but related-- quantity of interest is the exponent 
$\chi =\chi (d, \b)
\in [0,1/2]$
for the fluctuations of the normalizing constant, i.e. such that
$$
\ln Z_n - a_n  {\rm \ is \ of\ order \ } n^\chi {\rm \ for \ some \ 
constant \ } a_n 
$$
as $n \to \8$. 
These exponents $\xi, \chi$ depend also on 
the distribution of the environment $\h$. 

\medskip 

The ground state  of the model, defined as
the limit when $\b \to \8$, is the so-called oriented 
last passage percolation model. For the ground state  
it is believed that the exponents
$\xi(d, \8), \chi(d, \8)$ are universal, more precisely that they have the same
value for all distributions of  $\h$. 
Recently, Johansson, together with Baik and Deift, rigourously calculated the
values of these exponents in   dimension $d=1$ 
 and for  specific  distributions for $\h$. 
More precisely, in dimension $d=1$ 
 and for exponential and geometric 
 distributions, it is proven in \cite{Joh00}
that $\chi (1, \8)=1/3$, together with the Tracy-Widom law for 
limit fluctuations.
Also, for a one-dimensional Poissonized model, 
$\chi (1, \8)=1/3$ is obtained in \cite{BaDeJo99}
together with the Tracy-Widom  limit, though
$\xi(1, \8)=2/3$ is proved in \cite{Joh00b}:
the path is superdiffusive, in contrast with the underlying 
simple random walk
which is diffusive (corresponding to $\xi=1/2$).
\medskip

A number of predictions, conjectures and numerical estimates
can be found in the physical literature  \cite{KrSp91},
on the values on such exponents, and relations between them. In particular,
for all $\beta \in (0, \8]$, the scaling relation
$$
\chi=2 \xi-1\;
$$ 
is believed to hold in complete generality. This relation can be derived 
at a heuristic level
as a scaling in the Kardar-Parisi-Zhang equation \cite{KaPaZh86}, 
which status is, unfortunately, not clear
at a mathematical level. Instead,  partial results
 have been obtained rigourously  in specific situations
\cite{LiNePi96}, \cite{Piz97}, \cite{Mej04}, \cite{CaHu04}.
%%%%%%In fact, very little is known rigourously, especially for $d \geq 2$.
In fact, much is still open,  especially for $d \geq 2$.

\medskip

Bolthausen \cite{Bol89} placed the polymer model 
in the framework of martingales,
and noticed that the almost-sure limit of the rescaled partition function
is subject to  a dichotomy:
\bdnl{dych}
\lim_n \frac{Z_n}{ {Q}Z_n} \quad \left\{ 
\begin{array}{cc}
  >0 & \mbox{$Q$-a.s.}\\
  {\rm or} &\\
  =0 & \mbox{$Q$-a.s.}
\end{array}
\right.
\edn

A natural manner for measuring the disorder due to the random environment,
is to call {\bf weak disorder} the first case, and {\bf strong disorder}
the second one. Note that  weak disorder can be defined as the region 
where $\chi=0$ and $a_n=n \lm (\b )$.
The terminology is  justified by observing that 
the former case happens in large 
enough dimension for small $\b$ (including $\b
=0$) and the latter case 
for large $\b$ and general unbounded environment. More precisely,
a series of papers
\cite{ImSp88}, \cite{Bol89}, \cite{AlZh96}, \cite{SoZh96}
lead to the following.

\medskip

{\bf Theorem A} \quad {\it
%%\begin{theorem} \label{th1}(\cite{ImSp88} and \cite{Bol89}) 
Assume $d \geq 3$ and $\b$ small enough so that
\bdnl{numero}
%\frac{{Q}(\exp 2\b \h(t,x))}{({Q}\exp \b \h(t,x))^2}
%\times 
P( \exists n>0: \om_n=0)< \exp \left\{ -[\lm (2\b )-2\lm (\b) ] \right\}\;.
\edn
Then, weak disorder holds and, for almost every  
realization of the environment, the rescaled path:
\bdnl{om^(n)Int}
\om^{(n)}=\lef( \w_{nt}/\sqrt{n} \rig )_{t \ge 0},
\edn 
converges in law to 
the Brownian motion with diffusion matrix $d^{-1}I_d$.

%t \mapsto \hat\om^{(n)}_t:=n^{-1/2} \om_{nt}$, 
%%\end{theorem}
}
\medskip

This result was much a surprise for both mathematics and physics
communities who
did not expect that diffusivity could take place!

The second moment method was used to derive the theorem. The assumption 
on $\b$ is equivalent to the martingale $Z_n/QZ_n$ 
being bounded in $L^2$, and it is far 
from being necessary: A weaker quantitative condition
for weak disorder is obtained in \cite{Bir04} using size-biasing.   
Fifteen years were  necessary to improve on it:
The next result is a  criterion for weak disorder, where the 
critical quantity is
$$
  I_n= \mu_{n-1}^{\otimes 2} (\om_n=\tilde \om_n)\;,
$$
i.e., the probability for two polymers $\om$ and $\tilde \om$
independently sampled from the polymer measure in the same environment, 
to meet at time $n$.  

\medskip
{\bf Theorem B} \quad {\it
%%%\begin{theorem} \label{th2} 
(\cite{CaHu02} for the Gaussian case, \cite{CSY03}
for the general case). For non-zero $\b$ it holds

$$\Big\{ \lim_n \left(Z_n / {Q}Z_n\right) = 0\Big\} = \Big\{ \sum_n I_n = 
\8 \Big\}\quad
\mbox{$Q$-a.s.}
$$  
%%\end{theorem}
}
\medskip

The result is obtained by writing the semi-martingale decomposition of 
$\ln  Z_n / {Q}Z_n$, and studying separately the terms.
The above criterion is a
refined (conditional) second moment condition, and the criterion 
can also be used to obtain quantitative information on the polymer 
measure itself, on its concentration and localization
\cite{CSY03} in the strong disorder regime. 

In the present paper, we first establish the monotonicity in $\b$ 
concerning the dichotomy (\ref{dych}):
%%%%%%%%%%%%%
\Theorem{beta_cIntro}
%%%%%%%%%%%
There exists a critical value $\b_{\rm c}=\b_{\rm c}(d) \in [0,\8]$ with
\bdnn
\b_{\rm c}=0, 
& & \mbox{for $d=1,2$,} \\
0 < \b_{\rm c} \le \8 
& & \mbox{for $d \ge 3$}
\ednn
such that the weak disorder holds if $\b \in \{ 0\}\cup (0,\b_{\rm c})$ and 
the strong disorder holds if $\b >\b_{\rm c}$.
%%%%%%%%%%
\end{theorem}
%%%%%%%%%%%%
We also prove monotonicity for the Lyapunov exponent, see \Thm{beta_c}.
This result implies the absence of reentrant phase transition in the phase 
diagram of the model. The theorem follows from a correlation inequality
\cite{FKG71},
a natural ingredient which, however, appears here for the first time 
 (as far as we know)
in the field of directed polymers. 
\medskip

%Here, we will focus on the much simpler regime of weak disorder.
Now, we will focus on the regime of weak disorder.
There, it is natural to expect
that diffusive behavior takes place in the whole
weak disorder region, not only under the stronger assumption 
(\ref{numero}). Our main result is indeed:

\begin{theorem} \label{th3}
Assume $d \geq 3$ and weak disorder. Then, for all bounded
continuous functions $F$ on the path space, 
$$
\lim_n \mu_n[ F(\om^{(n)})]= {\bf E} F(B)
$$
in probability, where $\om^{(n)}$ 
is the rescaled path defined by (\ref{om^(n)Int}) and 
 $B$ is the Brownian motion  with diffusion matrix
 $d^{-1}I_d$. In particular, this holds for all $\b \in [0,\b_{\rm c})$.
\end{theorem}

Incidently, the statement shows that the scaling relation between exponents
does hold in the full weak disorder region, with $\xi=1/2$ and $\chi=0$. \\

In this paper, we also consider the fluctuations of extensive
thermodynamic 
quantities other than the partition function: we show that these are 
typically of order 1 -- like $\ln Z_n$ itself --, but we can prove this 
result only in part of the weak disorder region:

\begin{theorem} \label{th4}
Assume $d \geq 3$ and (\ref{numero}). Then, the 
energy averaged over the path
$$
 \mu_n[H_n] - n\lm'(\b ) \; \; 
%\frac{Q[\h(1,0)e^{\b \h(1,0)}]}{Q[e^{\b \h(1,0)}]} \quad 
\mbox{converges $Q$-a.s. to a finite random variable}.
$$
as $n \to \8$. A similar result holds for the entropy of $\mu_n$ 
with respect to $P$, see (\ref{eq:flentropy}).
\end{theorem}

In the proof of theorem \ref{th3} 
we introduce an infinite time horizon measure on the path space
which is a natural limit of the sequence $\mu_n$. This measure
is a time inhomogeneous Markov chain which depends on the environment.
We cannot prove the central limit theorem for this Markov chain
directly, but we need to average over the environment. In order to
prove convergence in probability with respect to the environment,
we use again a second moment method by introducing a second independent 
copy of the polymer before performing this average. All through, 
we  use the convergence of the series $\sum I_n$ as a main technical 
quantitative ingredient.

To prove Theorem \ref{th4} we use analytic functions arguments.
The crucial estimate is a bound on the second moment of
some complex random variable, this explains why we do assume  (\ref{numero}).
It is well known  that analytic
martingales are powerful tools to study disordered systems 
(e.g., section 5  of \cite{CN95}) in the regime of bounded second moment.

Our paper is organized as follows. After recalling some notations
and basic facts, we prove the existence of the critical temperature, 
together with characterization of the weak disorder phase 
that we will use further on
(section \ref{sec-CWD}).
We then introduce the Markov chain depending on the environment
in section \ref{sec-WD}. 
Section \ref{sec-CLT} deals with Gaussian behavior of the polymer,
and section \ref{sec-ana} with limits of energy and entropy. 
In the last section, we
illustrate the results in the case of Bernoulli environment,
emphasizing their relations with (last passage) oriented percolation.

\SSC{Notations and known facts}

Let 
$$
\bz_n(\om, \b)=\exp\{ \b H_n(\om) - n\lm(\b)\}
$$
Then, for all $\b$,
$$
W_n=Z_n \exp\{-n \lm(\beta)\} = Z_n/Q[Z_n] = P[ \bz_n(\om, \b)]
$$
is a positive 
martingale with respect to the $\s$-fields $\cG_n=\s\{ \h(s,x), 
s \leq n, x \in \Z^d\}$. By the martingale convergence theorem, 
it follows that 
$$
\lim_{n \ra \8}W_n=W_\8 \quad \mbox{$Q$-a.s.}\;,
$$
where $W_\8$ is a non-negative random variable. It is easy to see that 
the event $\{W_\8 >0\}$ is in the tail $\s$-field of $\{\cG_n, n\geq 0 \}$,
hence it is trivial by Kolmogorov 0-1 law. This shows the dichotomy
weak disorder versus strong disorder in
(\ref{dych}), which reads, in our new notation, 
\bdnl{dych2}
Q\{ W_\8 >0\} \; = \left\{ 
\begin{array}{cl}
  1 & \Llra \mbox{ weak disorder,}\\
  {\rm  or} &\\
  0 & \Llra \mbox{ strong disorder.}
\end{array}
\right.
\edn

It is well known \cite{CaHu02,CSY03}
that the weak disorder can happen if the transverse dimension
is large enough, i.e., $d \geq 3$. For $d\geq 3$,
\bdnl{eq:pid}
\pi_d :=P[ \exists n > 0: \om_n=0] \in (0,1)\;,
\edn
and (\ref{numero}) can be rephrased as
$$
\lambda (2\b)-2\lambda(\b) < - \ln \pi_d \quad
\Longrightarrow \quad W_\8>0\;\; \mbox{$Q$-a.s.}
$$
For $x \in\Z^d$, let $P^x$ be the law of the simple random walk 
in $\Z^d$ starting at $x$. If  $\theta_{n,x}$ denotes 
the shift operator given by  
$$\theta_{n,x} \h : (t,y) \mapsto 
\h (n+t,x+y)\;,$$  then  we have by  definition of $W_n$
$$
W_n \circ \theta_{0,x}=
P^x\left[ \bz_n \right]\;.
$$
By  definition of $W_n$ again,
 and by the simple Markov property,
we have also
\begin{equation}
  \label{eq:self}
W_n \circ \theta_{0,x}=
P^x\left[ \exp\{ \b \h(1,\om_1) - \lm(\b)\} W_{n-1} \circ \theta_{1,\om_1}
\right]
\end{equation}
and hence
\begin{equation}
  \label{eq:self8}
W_\8 \circ \theta_{0,x}=
P^x\left[ \exp\{ \b \h(1,\om_1) - \lm(\b)\} W_\8 \circ \theta_{1,\om_1}
\right]
\end{equation}
by taking the limit as $n \to \8$.
%%%%%%%%%%%%%%%%%%%%%%%%
\SSC{Characterizations of the  weak disorder phase and monotonicity}
\label{sec-CWD}
%%%%%%%%%%%%%%%%

We start by gathering some useful characterizations  of 
weak disorder, which should be compared to those in the case
where $\Z^d$ is replaced by a regular tree \cite[p.134]{kp}.
Before stating the next proposition, we make a remark.
For $\del \in (0,1)$, $(W_n^\del)$ is a uniformly integrable 
random variable. Therefore, 
%%the following limit exists always:
\bdnl{QW_8^del}
\lim_{n \ra \8}Q[W_n^\del]=Q[W_\8^\del].
\edn
%%%%%%%%%%%%%%%%%%%%
\Proposition{WDUI}
%%%%%%%%%%%%%%%%%%%%%%%%%%
 The following statements are equivalent for any 
$\del \in (0,1)$.
\bds
\item[(a1)] The martingale $W_n$ is uniformly integrable.
\item[(a2)] The martingale $W_n$ is $L^1$-convergent.
\item[(b1)] Weak disorder holds, i.e., $W_\8 >0$, $Q$-a.s.
\item[(b2)] The limit (\ref{QW_8^del}) is positive.
\item[(c1)] 
There exists a process
$$
(X_n,e_n)=((X_{n, x})_{x \in \zd}, (e_{n, x})_{x \in \zd}), \; \; n \in \N
$$
with values in $(\R^{\zd})^2$ such that
\bdmn
& & (e_n)_{n \in \N}\st{\rm law}{=}
\lef(\exp \{\b \h(n,\cdot )-\lm (\b)\}\rig)_{n \in \N },
\label{Law(e)=Law(eta)}\\
& & \mbox{For all $(n,x) \in \N \times \zd$, $Q[X_{n,x}]=1$,}
\label{Q[X_(nx)]=1}\\ 
& & \mbox{For all $(n,x) \in \N \times \zd$,  
${\dps X_{n, x}=P^x[e_{n+1,\w_1}X_{n+1, \w_1}] }$},
\label{X_n=P[e_(n+1)X_(n+1)]}\\
& & \mbox{For all $n\in \N $, $X_n$ is independent of $e_1,\ldots,e_n$.}
\label{X_nINDe_1,...,e_n}
\edmn
\item[(c2)] 
There exists 
%%an $\R^{\zd}$-valued random variable
  a non-negative random field 
$X=(X_x)_{x \in \zd}$
  on ${\zd}$
such that $Q[X_x]=1$ for all $x \in \zd$ and such that 
$$
X \st{\rm law}{=}
\lef(P^x[e_{\w_1}X_{\w_1}]\rig)_{x \in \zd}
$$
holds for any 
$\R^{\zd}$-valued random variable 
$e=(e_x)_{x \in \zd}$, independent of $X$, and 
$$
e \st{\rm law}{=}\exp (\b \h(1,\cdot )-\lm (\b)).
$$
\eds
%%%%%%%%%%%%%%%%%
\end{proposition}
%%%%%%%%%%%%%%%%%%
 \Remark{WDUI}
%%%%%%%%%%%%%%%%
Statements (a-1,2), (b-1,2) are natural. We will see in sections
\ref{sec-WD} and \ref{sec-CLT}, that 
(c1) is actually an important feature of the weak disorder
phase, allowing us to construct  the Markov chain $\mu$ in
(\ref{defmu}). 
The somewhat similar condition (c2) is  in the flavor of 
``condition ($\gm$)'' in \cite[Th\'eor\`eme 1 ]{kp}.
 \hfill $\Box$
\end{remark}

Proof of \Prop{WDUI}: 
(a1) $\Llra$ (a2): This follows from standard martingale convergence 
results \cite{durrett2}.\\
(b1) $\Llra$ (b2): This is obvious from the 
dichotomy (either $W_\8=0$, $Q$-a.s., or $W_\8>0$, $Q$-a.s.).\\
(a2) $\Lra$ (b1): The $L^1$-convergence implies $Q[W_\8]=1$, and hence (b1) 
by the dichotomy.\\
%%%To prove that (b1) implies (a1), we will go 
%%through the following condition (c1):
%%%%%%% \bds (c1) \eds
(b1) $\Lra$ (c1): 
Set 
\bdnl{X_ne_n}
X_{n,x}=W_\8 \circ \theta_{n,x}/Q[W_\8]\;, \quad 
e_{n,x}=\exp \{\b \h(n,x)-\lm (\b)\}.
\edn
We then have (\ref{Law(e)=Law(eta)}), (\ref{Q[X_(nx)]=1}) 
and (\ref{X_nINDe_1,...,e_n}). 
Moreover, 
we obtain (\ref{X_n=P[e_(n+1)X_(n+1)]}) by (\ref{eq:self8}). \\
(c1)$\Lra$ (a1): 
We will prove the uniform integrability by 
showing that 
\bdnl{W_n=Q[X|G]}
(Q[X_{0,0}|\tl{\cG}_n])_{n \ge 1}=(W_n)_{n \ge 1}
\edn
where $\tl{\cG}_n =\s [e_1,..,e_n]$.
Iterating (\ref{X_n=P[e_(n+1)X_(n+1)]}), we see from Markov property that
$$
X_{0,0}=P[e_{1,\w_1}\ldots e_{n,\w_n}X_{n,\w_n}].
$$
Taking the $Q$-expectation conditionally 
on $\tl{\cG}_n$, and 
observing (\ref{Q[X_(nx)]=1}) and (\ref{X_nINDe_1,...,e_n}), 
we arrive at
$$
Q[X_{0,0}|\tl{\cG}_n]=P^x[e_{1,\w_1}\cdots e_{n,\w_n}Q[X_{n,\w_n}]]
=P^x[e_{1,\w_1}\cdots e_{n,\w_n}],
$$
which proves (\ref{W_n=Q[X|G]}).\\

%%%%%%%%%%%%%%%%%%%%%%%%%%%%%%%%
%% \Remark{WDUI}
%%%%%%%%%%%%%
%% \bds (c2)  \eds

%% Proof of the equivalence:
(b1)$\Lra$ (c2): 
Define $X_n=(X_{n,x})_{x \in \zd}$ 
and $e_n=(e_{n,x})_{x \in \zd}$ by (\ref{X_ne_n}).
We prove that $X_1$ is what we look for. 
Since $X_1$ is independent of $e_1$, we have 
\bdnn
\lef(P^x[e_{\w_1}X_{1,\w_1}]\rig)_{x \in \zd}
& \st{\rm law}{=}&
\lef(P^x[e_{1,\w_1}X_{1,\w_1}]\rig)_{x \in \zd}\\
&=&\lef(W_\8 \circ \tht_{0,x}/Q[W_\8]\rig)_{x \in \zd}\\
& \st{\rm law}{=}&X_1,
\ednn
where we have used (\ref{eq:self8}) on the second line.\\
(c2)$\Lra$ (c1):
Suppose that $\ov{e}_n=(\ov{e}_{n,x})_{x \in \zd}$ ($n \in \N$) 
are independent of $X$ and 
$$
(\ov{e}_n)_{n \in \N}\st{\rm law}{=}
\lef(\exp (\b \h(n,\cdot )-\lm (\b))\rig)_{n \in \N }.
$$
We define 
$\ov{X}_n=(\ov{X}_{n,x})_{x \in \zd}$ ($n \in \N$) recursively by 
$$
\ov{X}_0=X, \; \; \ov{X}_{n+1,x}=P^x[\ov{e}_{n,\w_1}\ov{X}_{n,\w_1}].
$$
By the construction, $(\ov{X}_n, \ov{e}_n)$, $n=0,1,2,\ldots$ is a stationary 
process. Hence, the sequence of laws 
$$
\rh_n (ds_0 \cdots ds_n)=
Q((\ov{X}_{n-j}, \ov{e}_{n-j}) \in ds_j,\; j=0,\ldots,n ), 
\; \; \; n \in \N
$$
is consistent. Therefore, by Kolmogorov's extension theorem, there is a 
process $(X_n, e_n)$, $n=0,1,2,\ldots$ such that 
$$
Q((X_{j}, e_{j}) \in ds_j,\; j=0,\ldots,n )
=Q((\ov{X}_{n-j}, \ov{e}_{n-j}) \in ds_j,\; j=0,\ldots,n ),
\; \; \; n \in \N
$$
Then, (\ref{Law(e)=Law(eta)}) and (\ref{Q[X_(nx)]=1}) are obvious, while 
the recursion for $\ov{X}_n$ implies (\ref{X_n=P[e_(n+1)X_(n+1)]}). 
Finally, we see (\ref{X_nINDe_1,...,e_n}) 
from the fact that $\ov{X}_0$ and $\ov{e}_0,\ldots, \ov{e}_{n-1}$ are 
independent. \hfill $\Box$
%%%%%%%%%%%%
%%% \end{remark}
%%%%%%%%%%%

We now turn to the monotonicity of the phase transition.
We define the Lyapunov exponent by
\bdnl{psi(beta)}
\psi (\b )=-\limn \frac{1}{n}Q[\ln W_n]=\lm (\b)-\limn \frac{1}{n}Q[\ln Z_n].
\edn
The limit exists by subadditivity \cite[Proposition 1.5]{CSY03}. 
We see from Jensen's inequality that $\psi (\b )$ is non-negative. 
Moreover, $\psi$ is continuous 
in $\b$, since $\limn \frac{1}{n}Q[\ln Z_n]$ is convex in $\b$.
%%%%%%%%
\Theorem{beta_c}
%%%%%%%
\bds
\item[(a)] 
There exists a critical value $\b_{\rm c}=\b_{\rm c}(d) \in [0,\8]$ with
\bdmn
\b_{\rm c}=0, 
& & \mbox{for $d=1,2$,} \label{d=1,2}\\
0 < \b_{\rm c} \le \8 
& & \mbox{for $d \ge 3$}\label{dge3}
\edmn
such that 
\bdnl{beta_c}
Q\{ W_\8 >0\} = \left\{ 
\begin{array}{cl}
  1 & \mbox{if $\b \in \{0 \} \cup (0,\b_{\rm c})$,}\\
  0 & \mbox{if $\b >\b_{\rm c}$.}
\end{array}
\right.
\edn
\item[(b)] The Lyapunov exponent $\psi (\b)$ is non-decreasing in 
$\b \in [0,\8)$. In particular, there exists $\b^{\psi}_{\rm c}=\b^{\psi}_{\rm c}(d)$ with
\bdn
\label{b<tlb}
\b_{\rm c} \le \b^{\psi}_{\rm c} \le \8, 
\edn
such that 
\bdnl{tlbeta_c}
\psi (\b) \left\{
\begin{array}{cc}
  =0 & \mbox{if $ \b \in \R \cap [0, \b^{\psi}_{\rm c}]$,}\\
  >0 & \mbox{if $\b \in \R \bsh [0, \b^{\psi}_{\rm c}]$}
\end{array}
\right.
\edn
\eds
%%%%%%%%%%%%%
\end{theorem}
%%%%%%%%%%
%%%%%%%
\Remark{Unibeta_c}
%%%%%%%%%%
It is natural to expect that 
$\b_{\rm c}=\b^{\psi}_{\rm c}$, i.e., the absence of the 
intermediate phase. However, this is an open problem at the moment,
as well as whether weak or strong disorder hold at the critical value 
$\b_{\rm c}$.
%%%%%%%%
\end{remark}
%%%%%%%%%%%%%%%%
\Thm{beta_c} is a consequence of the monotonicity described in 
part (b) of the following lemma.   
%%%%%%%%%%%%
\Lemma{6/6beta}
%%%%%%%%%%%%
\bds
\item[(a)]
Assume that $\phi: (0,\8) \lra \R$ is  $\cC^1$ and 
that there are constants $C, p \in [1,\8)$ such that 
$$
|\phi'(u)| \le Cu^p+Cu^{-p}, \; \; \; \mbox{for all $u >0$.}
$$
Then, $\phi (W_n), \frac{\6 \phi (W_n)}{\6 \b} \in L^1(Q)$, 
$Q\phi (W_n)$ is $\cC^1$ in $\b \in \R$, and 
$$
 \frac{\6 }{\6 \b} Q \phi( W_n)= Q \frac{\6 }{\6 \b} \phi( W_n).
$$
\item[(b)]
Suppose in addition that $\phi$ is concave on $(0,\8)$. 
Then,
\bdnl{6/6beta}
Q \frac{\6 }{\6 \b} \phi( W_n) \le 0\; \; \; \mbox{for $\b \ge0$.}
\edn 
%%%%%%%%
% Q \frac{\6 }{\6 \b} \phi( W_n) 
%\lef.\barray{ll} \ge 0 & {\rm for } \; \b \geq 0\\  
%\le 0 & {\rm for } \; \b \leq 0 \earray \rig.
%%%%%%%%%%%%%%
\eds
%%%%%%%%%%%
\end{lemma}
%%%%%%%%%%%%%%%%%%%%%%%%
Proof: 
(a): Let $I=[0, \b_1]$ ($0 <\b_1<\8$)
and 
$$
X_n=\frac{\6 W_n}{\6 \b}=P[(H_n-n \lm') \bz_{n}].
$$
We first check that, for all $n$, 
\bds
\item[($\ast$1)] \hspace{1cm}
${\dps \sup_{\b \in I}W_n,\; \sup_{\b \in I}W_n^{-1},\; 
\sup_{\b \in I}|X_n| \in L^p(Q)}$ for all $p \in [1,\8)$, 
\eds 
and thereby that
\bds
\item[($\ast$2)] \hspace{1cm}
${\dps \sup_{\b \in I}\lef| \frac{\6 \phi(W_n)}{\6 \b} \rig|\in L^1(Q)}$.
\eds
For ($\ast$1), we have 
\bdnn
W_n^{-p} 
& \le & P[\bz_{n}]^{-p} \\
& \le & P[\bz_{n}^{-p}] \\
& \le & e^{pn\lm }P\exp 
\lef( p\b\sum_{1 \le s \le n} |\h(s,\om_s)|\rig).
\ednn
The property ($\ast$1) claimed for $W_n^{-1}$ is obvious from the 
above expression. $W_n^p$ and $|X_n|^p$ are bounded similarly.

The claim ($\ast$2) follows from ($\ast$1) and from
$$
\lef| \frac{\6 \phi(W_n)}{\6 \b} \rig|
=|\phi'(W_n)X_n| \le (CW_n^p+CW_n^{-p})|X_n|.
$$
It is now, easy to conclude part (a) of the lemma. Since 
$\phi(W_n)$ is $\cC^1$ in $\b \in \R$, we have 
$$
\phi(W_n(\b_1))=\phi (1)+\int^{\b_1}_0\frac{\6 \phi(W_n)}{\6 \b}d\b
\; \; \; \mbox{for all $\b_1 \in \R$.}
$$
The properties claimed in part (a) of the lemma follow from 
this expression, ($\ast$1) and Fubini's theorem.\\
(b): 
%%%%%%%%%%%%%%%%%%%%%%%%%%%%5
% previous proof is transferred (after "end{document}")
%%%%%%%%%%%%%%%%%%%%%%%%%%%%
We have 
$$
Q \frac{\6 }{\6 \b} \phi( W_n) = Q[\phi'(W_n)X_n]
= P\lef[ Q[\phi'(W_n)(H_n-n \lm') \bz_{n} ]\rig]
$$
Now, for a fixed path $\w$, 
%%%%%%%%%
%$H_n$ is linear in $\{\h (i,\w_i) \}^n_{i=1}$.
%With this observation, it is easy to see from 
%the argument in \cite{FKG71} that 
%%%%%%%%%%%%%%%%%%%%%%%%
the probability measure $\bz_{n}dQ$ is product, and therefore satisfies 
the FKG inequality \cite[p.78]{Lig85}. 
The function
$H_n-n \lm'$ is increasing in $\h$, while
$\phi'(W_n)$ is decreasing since $\phi$ is concave. 
These imply
$$
Q[\phi'(W_n)(H_n-n \lm') \bz_{n} ] \le 
Q[\phi'(W_n) \bz_{n} ]Q[(H_n-n \lm') \bz_{n} ]=0,
$$
and hence (\ref{6/6beta}).
\hfill $\Box$

\vvs
Proof of \Thm{beta_c}:
(a): By applying \Lem{6/6beta} to $\phi (x)=x^\del$ ($0<\del<1$),
it follows that the limit (\ref{QW_8^del}) is non-increasing 
in $\b \in [0,\8)$. This, together with \Prop{WDUI}, implies 
the existence of the values $\b_{\rm c}$ with the property 
(\ref{beta_c}). We then see (\ref{d=1,2}) 
from \cite[Theorem 1.3(b)]{CSY03}, and (\ref{dge3}) from Theorem A 
in section \ref{Intro}.\\
(b): By applying \Lem{6/6beta} to $\phi (x)=\log x$,
it follows that the limit (\ref{psi(beta)}) is non-decreasing 
in $\b \in [0,\8)$. This, together with the continuity of $\psi$,  implies 
the existence of the values $\b^\psi_{\rm c}$ with the property 
(\ref{tlbeta_c}). We then see (\ref{b<tlb}) from the obvious 
fact that $\psi (\b)>0$ implies $W_\8 =0$, $Q$-a.s.
\hfill $\Box$ 
%%%%%%%%%%
\SSC{The weak disorder polymer measure and its long time behavior}
\label{sec-WD}
%%%%%%%%%%%%%%%%

As a general fact, the measure $\mu_n$ is  a
(time-inhomogeneous) Markov chain, with 
transition probabilities 
$$
\mu_n(\om_{i+1}=y | \om_{i}=x)= \frac{
\exp\{ \b \h(i+1,y)-\lm\} 
 W_{n-i-1} \circ \theta_{i+1,y}}{
 W_{n-i} \circ \theta_{i,x}}\;
P(\om_{1}=y | \om_{0}=x)
$$
for $0 \leq i <n$, and $\mu_n(\om_{i+1}=y | \om_{i}=x)=
P(\om_{1}=y | \om_{0}=x)$ for $i \geq n$. Indeed, one can check that,
for any path $x_{[0,m]}=(x_0, \ldots x_m)$ of length $m \leq n$,
\begin{equation}
  \label{eq:munmarkov}
  \mu_n(\om_{[0,m]}=x_{[0,m]})= \bz_m
 \frac{W_{n-m}
\circ \theta_{m,x_m}}{{W_n}} P(\om_{[0,m]}=x_{[0,m]})\;.
\end{equation}

In the weak disorder
regime, we denote by $\mu$ the 
(random, time-inhomogeneous) Markov chain starting at 0 with 
transition probabilities 
\bdnl{defmu}
\mu(\om_{i+1}=y | \om_{i}=x)= \frac{\exp\{ \b \h(i+1,y)-\lm\} 
 W_{\8} \circ \theta_{i+1,y}}{
 W_{\8} \circ \theta_{i,x}}\;
P(\om_{1}=y | \om_{0}=x)\;.
\edn
In other respects,  for $A \in \cF_\8$ the limit
$$
\mu_\8(A) := \lim_{n \to \8} \mu_n(A)\;.
$$
exists by martingale convergence theorem for both numerator
and denominator of $\mu_n(A)$. 

The problem is that, it is not clear if 
%%%
the previous limit defines,  for a.e. $\eta$,
a probability measure on $\cF_\8$. But 
%$\mu_\8$ is a.s. $\sigma$-additive. On the other hand,
the Markov chain $\mu$ does. In the next result we relate 
these two objects $\mu_\8, \mu$, and we show that the latter
yields a nice description of the limit, in a precise sense.
%%%%%%%%%%%%%%%%%%%%
 \Proposition{prop-abscont}
%%%%%%%%%%%%%%%%%%%%%%%%
  Assume weak disorder. Then, 
\bdnl{m=m8cyl}
\mu(A)=\mu_\8(A)\;\mbox{$Q$-a.s. for $A \in \cup_{n \ge 1}\cF_n$.}
\edn
As a result, 
\bdnl{m_n->m}
Q \lef\{ \limn \m_n =\m \; \mbox{weakly}\rig\}=1.
\edn
Moreover, 
\bdmn
Q \mu (A)&=&Q \mu_\8 (A)\;,\quad \forall A \in \cF_\8\;,  
\label{trpp}\\
P &\ll& Q \mu \;\ll \; P \; \; \mbox{on $\cF_\8$.}\label{P<<Qm<<P}
\edmn
 \end{proposition}
%%%%%%%%%%%%%%%%%%%%%%%%%%%%%

\vvs
To prove \Prop{prop-abscont}, the following simple observation is useful.
%%%%%%%%%%
\Lemma{Amn}
%%%%%%%%%%
Suppose $\{ A_{m,n}\}_{m,n \ge 1} \sub \cF_\8$ are such that 
${\dps \lim_{m \nearrow \8}\sup_{n}P(A_{m,n})=0}$. Then
$$
\lim_{m \nearrow \8}\sup_{n}Q\m_n (A_{m,n})=
\lim_{m \nearrow \8}\sup_{n}Q\m_\8 (A_{m,n})=0.
$$
%%%%%%%%%%
\end{lemma}
%%%%%%%%%%%
Proof: 
We prove that ${\dps \lim_m\sup_{n}Q\m_n (A_{m,n})=0}$, 
the proof of the other one being similar.
For $\del >0$, 
$$
Q\m_n (A_{m,n})   \le 
Q[\m_n (A_{m,n}){\bf 1}_{W_n \ge \del} ]
+Q[W_n \le \del] 
$$
We have 
\bdnn
\sup_{n}Q[\m_n (A_{m,n}){\bf 1}_{W_n \ge \del}]
&\le & \del^{-1}\sup_{n}Q[W_n\m_n (A_{m,n})] \\
& = & \del^{-1}\sup_{n}P(A_{m,n}),
\ednn
which vanishes as $m \nearrow \8$. 
On the other hand, since $W_n^{-1}$ converges $Q$-a.s., 
their distributions are tight:
$$
\lim_{\del \searrow 0}\sup_{n}Q[W_n \le \del] =0.
$$
These prove the lemma.
%%%%%%%%%%%%%%%
 \hfill $\Box$
%%%%%%%%%%%%%%%%%

\vvs
%%%%%%%%%%%%%%%%%%%%%%%%%%%%%%%%%%%
Proof of \Prop{prop-abscont}: 
%%%%%%%%%%%%%%%%%%%%%%%%%%%%%%%%%%%
The first statement (\ref{m=m8cyl}) follows from 
(\ref{eq:munmarkov}). The second statement (\ref{m_n->m}) 
follows from (\ref{m=m8cyl}) by noting that the set of 
continuous functions on $\W$ contains a dense countable set 
of cylindrical functions.

To see (\ref{trpp}), we note that the averaged limit $Q \mu_\8 (A)$ is a 
probability measure on $\cF_\8$. Indeed, it is clearly finitely additive
by definition, and we have also by \Lem{Amn}, 
${\dps \lim_m Q \mu_\8 (A_m)=0}$ for 
any sequence $(A_m)_m$ in $\cF_\8$ which decreases to $\emptyset$.
Therefore, we have (\ref{trpp}) 
since the two probability measures $Q \mu$ and $Q \mu_\8 $ 
coincide on any $ \cF_n$. 

We see from \Lem{Amn} that $Q \mu \ll P$. 
To show the converse, assume
that  $Q \mu(A) \equiv Q \mu_\8(A)=0$. Then, $\mu_\8(A)=0$ a.s. and
$\mu_n(A) \to 0$ a.s. 
This implies that $W_n \mu_n(A)$ tends a.s. to 0 and,
combined with the uniform integrability of $(W_n)$,
it also implies that this sequence is itself uniformly 
integrable (recall $\mu_n(A) \leq 1$). Therefore,
$W_n \mu_n(A)$ tends  to 0 in $L^1(Q)$, that is,
$$
P(A)=Q[ W_n \mu_n(A)] \to 0\,, 
$$
which is the desired result. \hfill $\Box$

As a direct consequence, the polymer path inherits under $\mu$
the a.s. behavior of 
the simple random walk:
%%%%%%%%%%%%%%%%
\Remark{cor-ill}
%%%%%%%%%%%%%%%
 Assume weak disorder. Then, for $Q$-a.e. environment and $\mu$-a.e.
%% statement clear for  $\mu_
path,  
$$
\limsup_{n \to \8} \frac{\om_n}{\sqrt{2n \ln \ln n}}=1 \quad({\rm iterated\
 logarithm\ law})
$$
$$
\lim_n \frac{1}{\ln n}\sum_{ j \leq n}  \frac{1}{j} \delta_{\om_j/\sqrt j}
= {\cal N}(0, \frac{1}{d}{\rm Id}) 
%%%%%%\quad Q-a.s.
\quad({\rm a.s.\ central\ limit\ theorem})
$$
%%ALSO: Arc-sine law (cf berkes ??), as-Donsker,  AND
%%the continuous case with: continuity modulus and
%% set of zeros.. REFERENCES
\end{remark}

Regarding (\ref{m=m8cyl}), we have a more quantitative statement
concerning the variational norm $\| \nu-\nu' \|_{\cF_m}=
2 \sup\{ \nu(A)-\nu'(A); A \in {\cF_m}\}$. 
%%%%%%%%%%%%%%%%%%%%%%
\begin{proposition}
  \label{prop-muntomu}
%%%%%%%%%%%%%%%%%
 In the weak disorder case, 
$$
\lim_{k \to \8} \sup_m Q\left[ \|\mu_{m+k}-\mu\|_{\cF_m}\right]
=0
$$
\end{proposition}

\begin{remark} In particular, the central limit theorem for $\mu_n$
would follow from the one for $\mu$, but we could not prove the latter
directly.
\end{remark}

Proof of Proposition \ref{prop-muntomu}: 
We start to prove that
\begin{equation}
  \label{eq:lop1}
\sup_m
Q\left[ W_\8 \|\mu_{m+k}-\mu\|_{\cF_m}\right] \to 0\; , \quad k \to \8 \;.
\end{equation}
From (\ref{eq:munmarkov}) and the similar relation for $\mu$, for 
$m, k \geq 0$, it holds
\begin{eqnarray*}
 W_{\8}   \|\mu_{m+k}-\mu\|_{\cF_m}&=& W_{\8}
P\left[\bz_m \left|  \frac{W_{k}
\circ \theta_{m,\om_m}}{{W_{m+k}}} - \frac{W_{\8}
\circ \theta_{m,\om_m}}{{W_\8}} \right| \right]\\
&=& \frac{1}{W_{m+k}} P\left[\bz_m
\left| W_\8 W_{k}\circ \theta_{m,\om_m} -
W_{m+k} W_{\8} \circ \theta_{m,\om_m} \right| \right]\\
&\leq & \left| W_{\8}-W_{m+k} \right| +
 P\left[\bz_m
\left|  W_{k}\circ \theta_{m,\om_m} -
 W_{\8} \circ \theta_{m,\om_m}\right| \right]
\end{eqnarray*}
The   $Q$-expectation of the first term in the right-hand side
vanishes as $k \to \8$, though for the second one,
$$  Q\left(  P\left[\bz_m
\left|  W_{k}\circ \theta_{m,\om_m} -
 W_{\8} \circ \theta_{m,\om_m}\right| \right]\right)
= 
\qquad  \qquad  \qquad  \qquad \qquad  \qquad  \qquad  \qquad 
$$
\begin{eqnarray*}
 \qquad  \qquad &=& Q\left(  P\left[\bz_m
 Q( \left|  W_{k}\circ \theta_{m,\om_m} -
 W_{\8} \circ \theta_{m,\om_m}\right| | \cG_m)
\right]\right)\\
& =&  Q\left(  P\left[\bz_m \| W_{k}-  W_{\8}\|_{L^1(Q)}
\right]\right)\\
&=&\| W_{k}-  W_{\8}\|_{L^1(Q)} \longrightarrow_{k \to \8} 0\;.
\end{eqnarray*}
This proves (\ref{eq:lop1}). Now, it suffices to write
\begin{eqnarray*}
Q\left[ \|\mu_{m+k}-\mu\|_{\cF_m}\right]
&=&
Q\left[ \|\mu_{m+k}-\mu\|_{\cF_m} \left({\bf 1}_{W_\8 > \delta}+
{\bf 1}_{W_\8 \leq \delta}\right)
\right]\\
&\leq & 
\delta^{-1} 
Q\left[W_\8 \|\mu_{m+k}-\mu\|_{\cF_m}\right]+ 2 
Q\left[W_\8  \leq \delta \right]\;,
\end{eqnarray*}
and to optimize over positive $\delta$'s.
\hfill $\Box$

%%%%%%%%%%%%%%%%%%%%%%%%
\SSC{Central limit theorems}
\label{sec-CLT}
%%%%%%%%%%%%%%%%
Let $(\IW, \cF^{\IW}, P^\IW)$ be the $d$-dimensional Wiener space:
$$
\IW=\{ w \in C([0,1] \ra \rd)\; ; \; w(0)=0 \; \}
$$ 
with the topology induced by the uniform norm 
$\| w \|=\sup_{0 \le t \le 1}|w_t|$, let $\cF^\IW$  be the Borel 
$\s$-field and $P^\IW$ the Wiener measure .
For $n =1,2,\ldots$, we define the diffusive rescaling 
$\w \mapsto \w^{(n)}$ ($\W \ra \IW$) by 
\bdnl{om^(n)}
\w^{(n)}_t=\w_{nt}/\sqrt{n},\; \; \; 0 \le t \le 1,
\edn
where $(\w_{t})_{t \in \R_+} \in \IW$ is the linear interpolation of 
$(\w_{n})_{n \in \Z_+} \in \W$. This section is devoted to the proof 
of

%%%%%%%%%%%%%%%%%%%%%%%%
%%%%%%%%%%%%%
\Theorem{th:weakclt}
%%%%%%%%%%%%%%
Assume $d \ge 3$ and weak disorder. 
Then, for all $F \in C_{\rm b}(\IW)$, 
\bdmn
\limn \mu_n\left[ F(\w^{(n)})\right]
&=&P^\IW[F(w/\sqrt{d})], \label{weakclt1} \\
\limn \mu \left[ F(\w^{(n)})\right]
&=&P^\IW[F(w/\sqrt{d})], \label{weakclt2}
\edmn
in $Q$-probability.
In particular, these hold for all $\b \in [0,\b_{\rm c})$.
%%%%%%%%%%%%%%%
\end{theorem}
%%%%%%%%%%%%%%%%%%%%%%%%%%%%%%%%%%%%%%%%%%%%%%%%%%%%%%%%%%%%
\Remark{weakclt}
%%%%%%%%%%%
Since $F$ is bounded, 
the convergence in $Q$-probability claimed for (\ref{weakclt1}) and 
(\ref{weakclt2}) is equivalent to 
$L^p(Q)$-convergence for any finite $p$.
%%%%%%%%%%
\end{remark}
%%%%%%%%%%

\vs 
As a first step we start with the following weaker 
%%version of \Thm{th:weakclt}
statement, whose proof 
is also much simpler:
%%%%%%%%%%
\Proposition{ACLT}
%%%%%%%%%
Assume that weak disorder holds. Then,
\bdmn
\limn Q\m_n (\w^{(n)} \in \cdot )&=&P^\IW( w/\sqrt{d} \in \cdot ), 
\; \; \mbox{weakly.}\label{ACLT1}\\
\limn Q\m (\w^{(n)} \in \cdot )&=&P^\IW( w/\sqrt{d} \in \cdot ), 
\; \; \mbox{weakly.}\label{ACLT}
\edmn
%%%%%%%%%%
\end{proposition}
%%%%%%%%%%%
%%%%%%%%
\Remark{ACLT}
%%%%%%%
\begin{enumerate}
\item[(i)]
As can be seen from the proof below, (\ref{ACLT}) is true 
for any probability measure $R$ with $R \ll P$ instead of 
$Q\m$. 
\item [(ii)]
Of course, it is unnecessary to state and prove \Prop{ACLT} separately. 
However, the role of \Lem{lem:m2exists} below is made clearer in this way. 
\end{enumerate}
%%%%%%%%%
\end{remark}
%%%%%%%%%%%%
Proof: 
We write $\ov{F}(w)=F(w)-P^\IW[F(\cdot /\sqrt{d})]$ for $F \in C_{\rm b}(\IW)$.
We introduce the set $BL (\IW)$ of bounded Lipschitz 
functional on $\IW$ by 
$$ 
BL(\IW)=\{ F:\IW \ra \R \; ; \; \| F \|_{BL}\equiv \| F \|+\| F \|_L <\8\; \},
$$
where $\| F \|=\sup_{w \in \IW}| F(w )|$ and 
$$
\| F \|_L=\sup \lef\{ \frac{F(w)-F(\tl{w})}{\| w-\tl{w}\|}\; ; \; 
(w, \tl{w}) \in \IW \times \IW, \; w \neq \tl{w}\rig\}.
$$
{\it Step 1}: proof of (\ref{ACLT}).
As is well known, (\ref{ACLT}) is equivalent to that
\bdnl{ACLT2}
\limn Q\m [\ov{F}(\w^{(n)})]=0, \; \; \mbox{for all $F \in BL(\IW)$,}
\edn 
e.g., \cite[page 310, Theorem 11.3.3]{Dud89}. 
To show (\ref{ACLT2}), we make use of an almost 
sure central limit theorem for the simple random walk
in the following form. 
If $\{N_k\}_{k \ge 1} \sub \Z_+$ is an increasing sequence such that 
$\inf_{k \ge 1}N_{k+1}/N_{k}>1$, then for any fixed $F \in BL(\IW)$, 
\bdnl{11pripri}
\limn \frac{1}{n}\sum_{1 \le k \le n}\ov{F}(\w^{(N_k)})=0, 
\; \; \mbox{$P$-a.s.}
\edn 
This follows from the argument in \cite[pages 98 --100]{wa00}. 
Now, for any convergent subsequence of 
$a_n=Q\m [\ov{F}(\w^{(n)})]$, we can find a further subsequence 
$a_{N_{k}}$ with $\inf_{k \ge 1}N_{k+1}/N_{k}>1$. 
The point is that, by (\ref{P<<Qm<<P}),  
(\ref{11pripri}) holds with ``$P$-a.s.'' replaced by 
``$Q\mu $-a.s.'' Thus, by integrating, we obtain that
$$
\limn \frac{1}{n}\sum_{1 \le k \le n}a_{N_{k}}=0.
$$
Therefore, we necessarily have  (\ref{ACLT2}). 

{\it Step 2}: 
Now, we want to move from $\mu$ to $\mu_n$ in order to get (\ref{ACLT1}).
As before, we need only to prove that
\bdnl{ACLT3}
\limn Q\m_n [\ov{F}(\w^{(n)})]=0, \; \; \mbox{for all $F \in BL(\IW)$,}
\edn 
%%We assume for the moment that $F$ is, in addition, uniformly continuous. 
For $0 \leq k \leq n$, we write
\begin{eqnarray*}
  Q \mu_n \lef|\lef[\ov{F}(\w^{(n)}) \rig]\rig| 
&\leq& 
  Q\mu_n \lef| F(\w^{(n)})-F(\w^{(n-k)})\rig| \\
& +& \|\ov{F}\|  \sup_m Q\left[ \|\mu_{m+k}-\mu\|_{\cF_{m}}\right] \\
& + &  
  Q
 \mu \left\vert \left[ \ov{F}(\w^{(n-k)}) \right] \right\vert \;.
\end{eqnarray*}

As $n \to \8$ and for fixed $k$, the first and the last 
bounds vanish. In fact, we apply (\ref{ACLT}) to see that  the 
last bound vanishes. For the first one, we note that 
$F$ is uniformly continuous and that
$$
\sup_{\w \in \W}\max_{0 \le t \le 1}|\w^{(n)}_t
-\w^{(n-k)}_t|=O(k/\sqrt{n}) .
$$
Finally, letting 
$k \to \8$, the middle bound vanishes due to (\ref{prop-muntomu}).
This proves (\ref{ACLT1}).
%%%%%%%%% I cut the useless part
\hfill $\Box$

\vvs
The following lemma is a key to prove \Thm{th:weakclt}. 
%%%%%%%%%%%%%%%%%
\Lemma{lem:m2exists}
%%%%%%%%%%%%%%%%%%%%%
For all $B \in \cF_\8^{\otimes 2}$, 
the following limit exists a.s. in the weak disorder
region:
\bdnl{md(B)=}
\md(B)= \lim_{n \to \8} \mu_n^{\otimes 2}(B)\;.
\edn
Moreover,
\bdmn
\md(B) & = & \mu^{\otimes 2}(B)
\;,\quad \forall B \in \cup_{n \ge 1}\cF_n^{\otimes 2}, \label{trpp1} \\
Q\md(B) &=& Q\left[\mu^{\otimes 2}(B)\right] \;, \quad \forall B \in 
\cF_\8^{\otimes 2}\;, \label{trpp2} \\
Q\mu^{\otimes 2} &\ll & P^{\otimes 2} \quad {\rm on\ } \cF_\8^{\otimes 2} \;.
\label{abscont2}
\edmn
%%%%%%%%%%
\end{lemma}
%%%%%%%%%% 
\Remark{m2exists}
%%%%%%%%%%%%%%
It is tempting to think of $\md$ as ``$\mu_\8^{\otimes 2}$'',
but since we do not know if $\mu_\8$ is a.s. $\sigma$-additive,
the notation is not appropriate.
%%%%%%%%
\end{remark}
%%%%%%%%%%%%
\medskip

Proof:
Recall from Theorem B 
in section \ref{Intro} (Theorem 2.1 in \cite{CSY03}), that
the random series  $\sum_n I_n$ either converges  almost surely or
diverges  almost surely, according to weak or strong disorder.
We therefore have that
\bdnl{eq:WD}
\sum_n I_n < \8 \quad  \mbox{\rm $Q$-a.s.},
\edn
which will be the crucial estimate in the present proof.
%%is a main technical ingredient in the proof.

We start by proving that the limit (\ref{md(B)=}) exists.  
For a sequence $(a_n )_{n \geq 0}$ (random or non-random), we set 
$\D a_n=a_n-a_{n-1}$ for $n \geq 1$. 
For  $B \in \cF_\8^{\otimes 2}$ fixed,  
$ X_n \st{\rm def.}{=}P^{\otimes 2}[ \bz_n(\om) \bz_n(\tilde \om)
{\bf 1}_B]$ is a submartingale. The proof is based on the Doob's 
decomposition of the process $X_n$. 
We start by writting
\bdnl{X=M+A}
X_n=P^{\otimes 2}(B)+ M_n+A_n\;,
\edn
with $M_n$ a martingale, $M_0= A_0=0$, and $A_n$ the increasing process 
defined by its increments
\bdmn
\D A_{n}&=& Q[\D X_n | {\cG_{n-1}}] \nn \\
&=& Q\left[
P^{\otimes 2}\left[ \bz_{n-1}(\om) \bz_{n-1}(\tilde \om){\bf 1}_B
\left\{ e(n, \om_n)e(n, \tilde \om_n)-1\right\}
\right] 
\Big| {\cG_{n-1}} \right]
\nn \\
&=& c P^{\otimes 2}\left[ \bz_{n-1}(\om) \bz_{n-1}(\tilde \om){\bf 1}_B
{\bf 1}_{\om_n=\tilde \om_n}
\right] \nn \\
&=&
c W_{n-1}^2 \mu_{n-1}^{\otimes 2}(B \cap \{\om_n=\tilde \om_n\})
\label{DAn=}\\
&\leq& c W_{n-1}^2 I_n \;,\label{DAn<}
\edmn
where $e(n, x)=\exp \{ \beta \h(n,x)-\lm(\b)\}$ and
the constant $c=\exp\{\lm(2 \b)-2 \lm(\b)\}-1$ is finite.
Hence the increasing process converges,
\bdnl{A8<}
A_n \nearrow A_\8 \leq c (\sup_k W_k)^2 \sum_k I_k < \8 
\; \; \; \mbox{$Q$-a.s}.
\edn
We prove that the martingale $(M_n)$ converges $Q$-a.s. by showing 
that 
\bdnl{<M><8}
\lan M \ran_\8 \le 4C(\sup_k W_k)^4 \sum_k I_k < \8\; \; \; \mbox{$Q$-a.s.}
\edn
with a constant $C=C(\b)$. 
Introducing 
$$
\vp_n (\w, \tl{\w})=e(n, \om_n)e(n, \tilde \om_n)-1
-c {\bf 1}_{\om_n=\tilde \om_n}, 
$$
we have by (\ref{X=M+A}), (\ref{DAn=}), 
\bdmn
\D M_{n}&=& 
P^{\otimes 2}\left[ \bz_{n-1} (\om)  \bz_{n-1} (\tilde \om) 
\vp_n (\w, \tl{\w}){\bf 1}_B\right] \nn \\
&=& 
W_{n-1}^2 \mu_{n-1}^{\otimes 2}[\vp_n (\w, \tl{\w}){\bf 1}_B], \label{DMn=}
\edmn
and hence 
\bdmn
\D \lan M \ran _{n}&=& 
Q \left[(\D M_{n})^2 \Big| {\cG_{n-1}}
\right] \nn \\
&=& 
W_{n-1}^4 \mu_{n-1}^{\otimes 4}[
Q[\vp_n (\w^1, \w^2)\vp_n (\w^3, \w^4)]{\bf 1}_{B \times B}],
\edmn
where $\w^1, \ldots, \w^4$ are independent copies of the path $\w$.
Note that $Q[\vp_n (\w, \tl{\w})]=0$ and that 
%%$C\st{\rm def.}{=}\sup_n Q[\vp_n (\w, \tl{\w})^2]$ is a constant
%%depending
%%%%%%%%%%%%%%%%%% FC changed
$C\st{\rm def.}{=}\sup\{ Q[\vp_n (\w, \tl{\w})^2]; n, \w, \tl{\w}\}$ 
is a finite constant depending
only on $\b$. We see from these and Schwarz inequality that
\bdmn
\D \lan M \ran _{n} & \le &
C\sum_{\scriptstyle i=1,2 \atop \scriptstyle j=3,4}
W_{n-1}^4 \mu_{n-1}^{\otimes 4}\left[\{\w_n^i=\w_n^j\} 
\cap (B \times B)\right] \label{D<M><} \\
& \le & 4CW_{n-1}^4I_n,\label{D<M><2}
\edmn
leading to (\ref{<M><8}).
This proves that  $X_n$,
as well as  $\mu_n^{\otimes 2}(B)=W_n^{-2}X_n$, 
converges $Q$-a.s. 
\medskip

As for (\ref{trpp1}), it follows from (\ref{eq:munmarkov})
directly that $\md$ and $\mu^{\otimes 2}$ coincide on cylindric 
events. 
%%%%%%%%%%%%%%%
%Here comes the key point of the argument.
%%%%%%%%%

As in (\ref{trpp}), the claims (\ref{trpp2}) and (\ref{abscont2}) 
boil down to proving 
that 
$$
\lim_{m \to \8} Q \md (B_m) 
=0 \;,
$$
for any $\{B_m \} \sub \cF^{\otimes 2}$ 
with $\lim_{m \to \8} P^{\otimes 2} (B_m) =0$. 
It is enough to prove that
$$
\lim_{m \to \8}  \md (B_m) \equiv \lim_{m \to \8}  \lim_{n \to \8}
\mu_n^{\otimes 2} (B_m) =0 \quad \mbox{in $Q$-probability},
$$
and hence that
\bdnl{1e}
\lim_{m \to \8}\sup_n X^{(m)}_n=0 \quad \mbox{in $Q$-probability},
\edn
where $X^{(m)}_n=P^{\otimes 2}
[ \bz_n(\om) \bz_n(\tilde \om){\bf 1}_{B_m}]$.
Let 
$$
X^{(m)}_n=
P^{\otimes 2}(B_m) + M_n^{(m)}+ A_n^{(m)} 
$$
be the submartingale decomposition as (\ref{X=M+A}).
Of course, $\lim_{m \nearrow \8} P^{\otimes 2}[B_m]=0$. 
Observe that, similar to (\ref{A8<}), it follows from
(\ref{DAn=}) that
$$
A_n^{(m)} 
 \leq c (\sup_k W_k)^2 P^{\otimes 2}(S{\bf 1}_{B_m})\;,$$
where 
$$S= \sum_{n\geq 1} \bz(\om, n-1)\bz(\tilde \om,n-1) 
{\bf 1}_{\om_n=\tilde \om_n}\;.$$
Now, the weak disorder assumption (\ref{eq:WD}) states that
this variable $S$ is $P^{\otimes 2}$-integrable for 
$Q$-almost every environment. Therefore,
\bdnl{A(m)->0}
\lim_{m \nearrow \8} A_\8^{(m)} =0,
\; \; \; \mbox{$Q$-a.s.}
\edn
For $M_n^{(m)}$, we see from (\ref{D<M><}),(\ref{D<M><2}) and 
the weak disorder assumption (\ref{eq:WD}) that 
\bdnl{<M(m)>->0}
\lim_{m \nearrow \8} \lan M^{(m)} \ran_\8 =0,
\; \; \; \mbox{$Q$-a.s.}
\edn
This implies that
\bdnl{M(m)->0}
\lim_{m \nearrow \8}\sup_n | M_n^{(m)}|=0 \; 
\mbox{in $Q$-probability} \;.
\edn
In fact, let $\t (\ell )=\inf \{ n \ge 0\; ; \; 
\lan M^{(m)} \ran_{n+1} >\ell \}$.
Then, 
$$
Q\{ \sup_n | M_n^{(m)}| \ge \e \} 
 \le  
Q\{ \lan M^{(m)} \ran_\8 >\ell \}
+Q\{ \sup_n | M_n^{(m)}| \ge \e, \;  \t (\ell )=\8\}.
$$
Clearly, the  first term on the right-hand-side vanishes 
as $m \nearrow \8$ and so does the second term as can be seen from 
the following application of Doob's inequality:
\bdmn
Q\{ \sup_n | M_n^{(m)}| \ge \e, \;  \t (\ell )=\8\}
& \le & Q\{ \sup_n | M^{(m)}_{n \wedge \t (\ell )}| \ge \e \} \nn \\
& \le & 4\e^{-2} Q[ \lan M^{(m)} \ran_{\t (\ell )}] \nn \\
& \le & 4\e^{-2} Q[ \lan M^{(m)} \ran_\8 \wedge \ell ] \nn
\edmn

By, (\ref{A(m)->0}) and (\ref{M(m)->0}), 
we conclude (\ref{1e}). \hfill $\Box$

\vs 
{\bf Proof of \Thm{th:weakclt}: } 
We write $\ov{F}(w)=F(w)-P^\IW[F(\cdot /\sqrt{d})]$ for $F \in C_{\rm b}(\IW)$.
We begin by proving (\ref{weakclt2}).
Repeating the same argument as in the step 1 of the proof of
Proposition \ref{Prop.ACLT},
%%By the same argument as in \Prop{ACLT}, 
but
using (\ref{abscont2}) 
instead of (\ref{P<<Qm<<P}), we obtain 
\bdnl{limG(ww)}
\limn Q\mu^{\otimes 2}[G(\w^{(n)},\tl{\w}^{(n)})]
=(P^\IW) ^{\otimes 2}[G(w/\sqrt{d},\tl{w}/\sqrt{d})]
\edn
for any $G \in C_{\rm b} (\IW \times \IW)$. Now, if we take 
$G(w,\tl{w})=\ov{F}(w)\ov{F}(\tl{w})$, 
then (\ref{limG(ww)}) reads
$$
\limn Q \lef[ \lef(\mu \left[ \ov{F}(\w^{(n)}) \rig] \rig)^2\rig]=0,
$$
which proves (\ref{weakclt2}). 

To obtain (\ref{weakclt1}) from (\ref{weakclt2}), 
we show that 
$$
\limn Q\lef | \m_n [\ov{F}(\w^{(n)})]\rig|=0
\; \; \; \mbox{for all $F \in C_{\rm b}(\IW)$.}
$$
This can be done 
by exactly the same approximation procedure as we used to deduce (\ref{ACLT1}) 
from (\ref{ACLT}), see  step 2 in the proof of
Proposition \ref{Prop.ACLT}.
%% where we proved that
%% $$ \limn Q \m_n \left[ \ov{F}(\w^{(n)})\right]=0
%% \; \; \; \mbox{for all $F \in C_{\rm b}(\IW)$}.$$ 
%% 
\hfill $\Box$

%%%%%%%%%%%%%%%%%%%%%%%%
\SSC{An analytic family of martingales}
\label{sec-ana}
%%%%%%%%%%%%%%%%

For $\b$ complex, $Q[ \exp \b \h(n,x) ]$ is well defined, but
we also want its logarithm to be holomorphic.
Let $U_0$ be the open set in the complex plane given by
$$
U_0=
{\rm \ connected \ component \ of \ } 0 {\rm \ in  \ }
\{ \b \in \C ;\; Q[ \exp \b \h(n,x) ] \notin  \R_-\}\;.
$$
Then, $U_0$ is a neighborhood of the real axis, and
$ \lm (\b)= \log Q[ \exp \b \h(n,x) ]$ is an analytic function on $U_0$.
Define, for $n \geq 0$ and $ \b \in U_0$,
\begin{equation}
  \label{defmartc}
  W_n(\b) = P\left[ \exp \left( \b \sum_{t=1}^{n} \h(t,\om_t)
    - n \lm(\b)\right) \right] \;.
\end{equation}
Then, for all $ \b \in U_0$, the sequence $(W_n(\b), n \geq 0)$
is a $(\cG_n)_n$-martingale with complex values, and for fixed $n$,
$W_n(\b)$ is an analytic function of $\b \in U_0$.

In view of the implication below (\ref{eq:pid}), we introduce for
%%Recall now the following sufficient condition for weak disorder
%%\cite{Bol89}. Assuming  
$d \geq 3$, the real subset 
\begin{equation}
  \label{eq:L2}
U_1=
\Big\{\b \in \R \;:\;  \lm(2\b)-2\lm(\b) < - \ln \pi_d
%%P[ \exists n > 0: \om_n=0]
\Big\}\;,
\end{equation}
which
is the set of $\b \in \R$ such that the martingale $(W_n)_n$ is $L^2$-bounded.
It is an open interval such that 
$0 \in U_1 \sub \{\b \in \R \; ;\;  W_\8(\b)>0\}$, $Q$-a.s.  

\begin{proposition} \label{prop:martc}
Assume  $d \geq 3$. 
%%the set 
%% \begin{equation}
%%  \label{eq:L2C}
%%U_2= \Big\{\b \in U_0 \;:\;   \lm(2 \Re \b)-2\Re\lm( 
%%\b) < - \ln P[ \exists n > 0: \om_n=0]\Big\}
%%\end{equation} 
Define $U_2$ as the connected component of the set 
$$
\Big\{\b \in U_0 \;:\;   \lm(2 \Re \b)-2\Re\lm( 
\b) < - \ln \pi_d
%%P[ \exists n > 0: \om_n=0]
\Big\}
$$
which contains the origin. Then, $U_2$ 
is a complex neighborhood of $U_1$, such that, as $ n \to \infty$,
$$
 W_n(\b) \to W_\8(\b)\;,\quad \mbox{$Q$-a.s.},
$$
where the convergence holds in the sense of analytic function.
In particular, the limit 
 $ W_\8(\b)$ is holomorphic in $U_2$, and $Q$-a.s.,
$$
\frac{d^k}{d\b^k} W_n(\b) \to \frac{d^k}{d\b^k} W_\8(\b) \;,
%%\quad k \geq 0
$$
uniformly on compacts of   $U_2$ ($ k \geq 0$).
\end{proposition}
%%%%%%%%%%%%%%%%%%%%%%%%%%%%%%%%%%%%%%%%%%%%%%%%%%%%%%%%%%%%
Proof of Proposition \ref{prop:martc}: From 
$\overline{(e^z)}=e^{\overline{z}}$ and $\overline{Q[f]}=
Q[\overline{f}]$, we see that 
$\overline{\lm (\b)}=\lm (\overline{\b})$, and that 
\begin{eqnarray}
Q\Big[ |W_n(\b)|^2\Big]&=& Q \Big[  P[
\exp\{ \b H_n(\om) -  n  \lm(\b)\}]
 P[\exp\{ \overline{\b} H_n(\tilde \om)-  n  \overline{\lm(\b)}\}]
\Big] \nonumber
\\\nonumber
&=&   P^{\otimes 2}\Big[ Q \Big[ 
\exp\{ \b H_n(\om) + \overline{\b} H_n(\tilde \om)- 2 n \Re \lm(\b)\}
\Big]\Big]\\\nonumber
&=&  P^{\otimes 2}\Big[ \exp\{  [\lm(2 \Re \b)-2\Re\lm( 
\b)]
\sum_{t=1}^n {\bf 1}_{\om_t=\tilde \om_t}
\}\Big]\\
& \nearrow& P^{\otimes 2}\Big[ \exp\{ [\lm(2 \Re \b)-2\Re\lm( 
\b)]
\sum_{t=1}^\8 
{\bf 1}_{\om_t=\tilde \om_t}
\}\Big] < \8
\label{matmata}
\end{eqnarray}
if $\b \in U_2$. 
%%For further purposes, let us note also that $U_2$
%%is connected, since for $\b \in U_2$ we have  $\Re \b \in U_2$
%%and the whole segment between these two points is in  $U_2$.

Now, let a point $\b \in U_2$, a radius $r>0$ such that the closed disk
$D(\b, r) \subset U_2$. Choosing $\rho >r$ such that $D(\b, \rho) 
\subset U_2$, we obtain by Cauchy's integral formula for
all $\b' \in D(\b, r)$,
$$
W_n(\b')=  \frac{1}{2i\pi}
\int_{\partial D(\b, \rho)} \frac{W_n(z)}{z-\b'} dz
  =  \int_0^1 
\frac{W_n(\b+\rho e^{2i\pi u})\rho e^{2i\pi u}}{(\b+\rho e^{2i\pi u})-\b' } 
du\;,
$$
hence
$$
X_n:=\sup\{ |W_n(\b')|; \b' \in  D(\b, r)\} \leq
\rho  \int_0^1  \frac{|W_n(\b+ \rho  e^{2i\pi u})|}{\rho-r} du
$$
Letting $C=(\rho /(\rho-r))^2$, we obtain by Schwarz inequality
  \begin{eqnarray*}
(Q[ X_n])^2 &\leq& C Q[ \int_0^1 |W_n(\b+ \rho  e^{2i\pi u})|^2 du]\\
    &\leq & C \sup\{  Q[ |W_n(\b'')|^2]; n \ge 1, \b'' \in  D(\b, \rho)\}\\
    & <& \8
  \end{eqnarray*}
in view of (\ref{matmata}). Notice now that $X_n$, a supremum of 
positive submartingales, 
 is itself  a positive submartingale. Since $\sup Q[ X_n]< \8$, 
$X_n$ converges Q-a.s. to a finite limit $X_\8$. Finally,
$\sup\{ |W_n(\b')|; \b' \in  D(\b, r)\}<\8$ a.s., and 
$W_n$ is uniformly bounded on compact subsets of $U_2$
on a set of environments of full probability. On this set, 
$(W_n, n \geq 0)$ is a normal sequence \cite{rudin}
which has a unique limit on the real axis: Since $U_2$ is
connected,
the sequence converges to some limit
$W_\8$, which is holomorphic on $U_2$, and the derivatives 
also converges to those of $W_\8$. \hfill $\Box$

Note that we do not know that $W_\8(\b) \neq 0$ for general
$\b \in U_2$, except for $\b \in U_1$ -- and of course for some complex 
neighborhood around $U_1$ --. We draw now some consequences 
for real $\b$'s. We write $ \mu_n= \mu_n^{\b}$ to recall the dependence
on the temperature.
\begin{theorem} \label{th:flenergy}
Assume $d \geq 3$. Then $W_\8$ and $\ln W_\8$ are analytic (real) 
function of $\b \in U_1$. Moreover, as $n \to \8$, 
\begin{equation}
  \label{eq:flenergy}
  \mu_n^{\b} [ H_n] - n \lm'(\beta) \to (\ln W_\8)'(\beta)\;,
\end{equation}
though for the entropy $h(\mu_n^{\b} | P)=\mu_n^{\b}[\ln (d\mu_n^{\b}/dP)]$,
\begin{equation}
  \label{eq:flentropy}
h(\mu_n^{\b} | P)-n[\beta  \lm'(\beta)- \lm(\beta)]
\to \beta(\ln W_\8)'(\beta)-\ln W_\8(\beta)\;,
\end{equation}
for all $\b \in U_1$.

On the other hand, for $Q$-a.e. environment, 
\begin{center}
  the law of $\frac{\displaystyle H_n-n\lm'(\b)}{\displaystyle \sqrt n}$
under $\mu_n$ converges to the Gaussian $\cN(0, \lm''(\b))$
\end{center}
where $\lm''(\b)>0$.
%%%
%%%In other respects, for all $u \in \R$, 
%%%\begin{equation}
%%%  \label{eq:flnorm}
%%%  \mu_n[ \exp \{ un^{-1/2}(H_n-n\lm'(\b))\}] \to \exp\{ -\lm''(\b) u^2/2\}
%%% \quad  $Q$-a.s.
%%%\end{equation}
%%%where $\lm''(\b)>0$.
\end{theorem}
{\bf Comment}: The {\it average energy} for the polymer measure,
$ \mu_n^{\b} [ H_n]$,
scales like the annealed one $ n \lm'(\beta)$, but it has 
fluctuations of order one in this part of the weak disorder region.
 The  entropy also has $\cO (1)$ fluctuations. On the other hand, 
the last result shows that,  due
to variations from a path to another,
the fluctuations of the  {\it energy under the polymer measure}
is normal and of order of magnitude $\cO(\sqrt n)$. \hfill $\Box$
\medskip
%%%%%%%%%%%%%%%%%%%%%%%%%%%%%%%%%%%%%%%%%%%%%

Proof of Theorem \ref{th:flenergy}:
We have the identities
\begin{eqnarray*}
   (\ln W_n)'(\beta) = \mu_n^{\b} [ H_n] - n \lm'(\beta),\\
   h(\mu_n^{\b} | P)=\beta \mu_n^{\b} [ H_n]-n \lm'(\beta)
-\ln W_n(\beta)\;.
\end{eqnarray*}
In view of  Proposition \ref{prop:martc},
 $(\ln W_n)'(\beta)=  ( W_n)'(\beta)/ W_n(\beta)$
converges a.s. to $ ( W_\8)'(\beta)/ W_\8(\beta)=(\ln W_\8)'(\beta)$
for  $\b \in U_1$,
which is the first result (\ref{eq:flenergy}). The second one 
(\ref{eq:flentropy}) follows 
easily. In order to prove the last one, we show the stronger statement
 that, for $Q$-a.e. environment,
$$
%%%\label{eq:flnorm}
 \mu_n\left[ \exp \{ 
\frac{\displaystyle u(H_n-n\lm'(\b))}{\displaystyle \sqrt n}
%un^{-1/2}(H_n-n\lm'(\b))
\}\right] \to \exp\{ \frac{\lm''(\b) u^2}{2}\}
$$
as $n \to \8$ for all $u \in \R$ and $\beta \in U_2$.
Write the left-hand side as
$$
\frac{W_n(\b + un^{-1/2})}{W_n(\b)}\times \exp  \Big\{n[\lm(\b+ un^{-1/2})-
\lm(\b)-un^{-1/2}\lm'(\b)]\Big\}
$$
Since $W_n \to W_\8$ locally-uniformly on $U_1$, and since $\lm$ is smooth,
the right-hand side converges $Q$-a.s. to 
$[W_\8(\b)/W_\8(\b)] \times\exp\{ \lm''(\b) u^2/2\}$ as $n \to \8$.
\hfill $\Box$

%%%%%%%%%%%%%%%%%%%%%%%%%%%%%%%%%%%%%%%%
%%\begin{example}{\bf Bernoulli environment}.
%%%%%%%%%%%%%%%%%%%%%%%%
\SSC{Bernoulli environment}
\label{sec-Bern}
%%%%%%%%%%%%%%%%
\rm Let $p\in (0,1)$. In this section, we focus on the Bernoulli case, where 
%% With $p\in (0,1)$, let 
$$
\h(t,x) =   
\left\{
    \begin{array}{c}
     0 \\ -1 
  \end{array}
\right. {\rm \ with\ }Q-{\rm probability \ }
\left\{
 \begin{array}{c}
     p \\ 1-p 
  \end{array}
\right. \;,
$$
%Q[\h(t,x)=0]=p=1-Q[\h(t,x)=-1]$, 
In this case,
$
\lm(\b)=\ln[p+(1-p)e^{-\b}]
$.

Consider also the  site, oriented Bernoulli percolation (see \cite{durrett1},
\cite{grimmett}), as follows:
Call a  site $(t,x)\in \N \times \Z^d$ open if $\h(t,x)=0$, and closed
if $\h(t,x)=-1$.
Write  $(n,x) \to^{\h} (k,z)$ if there exists an oriented open path
$((t,\om_t); n \leq t \leq k)$ from  $(n,x)$ to $ (k,z)$, i.e.,
some path $((t,\om_t); n \leq t \leq k)$ with nearest neighbors vertices 
$\om_t$ and $\om_{t+1}$ and $\h(t, \om_t)=0$ for all $t$, and   
$\om_n=x, \om_k=z$. Write $(n,x) \to^\h \8$ if there exists an infinite
oriented open path starting at $(n,x)$, and denote by 
$\cC$ the set of sites $(n,x)$ such that  $(n,x) \to^\h \8$ and
$\|x\|_1 \leq n$, $\|x\|_1 =n$ modulo 2. The set  $\cC$ is called  the 
infinite cluster. It is well known that 
there exists some  percolation
threshold  $\vec p_c(d) \in (0,1)$  such that 
for $p>\vec p_c(d)$ and $d \geq 1$, $\cC$ is $Q$-a.s. non empty,
and $\cC$ is $Q$-a.s.  empty for  $p<\vec p_c(d)$.
It is known 
(Theorem 2 in \cite{gh}), that 
$\cC$ is a.s. connected, in the sense that a.s. on the set  $\{(n,x) \to \8,
(m,y) \to \8\}$, there exists some $(k,z) \to \8$ such that both
 $(n,x) \to (k,z)$ and 
$(m,y) \to (k,z)$. Let  $H_n^*$ be
the maximum value of $H_n$ over all paths $\om$ starting from $(0,0)$.
In the last passage percolation problem, one is interested in 
the almost-sure limit 
$$
\tau=\limn-H_n^*/n,
$$ 
(called the time constant),
which exists and is constant by subadditivity \cite{durrett1},
\cite{grimmett}, and is non-negative. For directed polymers on the other hand, 
the a.s.-limit $\psi(\beta)=\psi(\beta,p)$ of $-(1/n)\ln W_n(\beta)$  
exists, is constant
 by  subadditivity and concentration \cite{CSY03}, and is non-negative.

We have a commutative diagram, with $\beta, n$ tending to $+\8$:
$$
\begin{array}{ccc}
-\frac{1}{n\b} \ln W_n(\b)  & \longrightarrow_\beta & -\frac{H_n^*}{n}\\
\downarrow_n &&\downarrow_n \\
\frac{\psi(\b)}{\b}  & \longrightarrow_\beta & \tau
\end{array}
$$
The proofs of the horizontal limits are easy, and left to the reader. 
We have $\tau=0$ for  $p>\vec p_c(d)$
by definition of the percolation
threshold,
and $\tau >0$ for  $p<\vec p_c(d)$
in view of the exponential tails of the cluster of the origin
\cite{menshikov}. Let us introduce another critical value,
\begin{equation}
  \label{eq:pcpsi}
  p_c^\psi= \inf\{ p \in [0,1]:  \psi(\beta; p)=0 \;,\forall
\b >0 \}
\end{equation}
and recall $\pi_d$ from (\ref{eq:pid}). We have
\begin{equation}
  \label{eq:pcpsi>}
  \pi_d \geq p_c^\psi \geq \vec p_c(d)\;.
\end{equation}
Indeed, it holds
$$
\frac{
\psi(\b)}{\b}=\frac{
\lm(\b)}{\b}-\lim_n\frac{1}{n\b}
\ln P[\exp \b H_n] 
%%\leq \frac{\lm(\b)}{\b}-\lim_n\frac{H_n^*}{n}+ \frac{d}{\b}\ln 2\;,
\geq \frac{\lm(\b)}{\b}+\tau\;,
$$
which
 becomes strictly positive in the limit $\b \to \8$ if $p<\vec p_c(d)$.
This proves the second inequality in our claim (\ref{eq:pcpsi>}). Now,
the first one follows from the observation that $[0, \8[ \subset U_1   $
holds if $p >\pi_d$, 
 (see example 2.1.1 in \cite{CSY04} for instance).

\medskip 

From now on, we assume that $d \geq 3$ and
$$p> \pi_d \;,
$$
which implies that weak disorder holds for all $\b \geq 0$
and also $\tau=0$.
There are a strong analogies between our limiting fluctuations  
from the previous section for the directed polymer model,
and the first passage time in oriented percolation. We now elaborate 
on these relations.

 We have the identities
 \begin{equation}
   \label{eq:887}
\lim_{\b \to +\8} \mu_n^\b[H_n]=\lim_{\b \to +\8}
\frac{1}{\b} \ln W_n= H_n^*\;,
 \end{equation}
 \begin{equation}
   \label{eq:888}
\lim_{n \to +\8} H_n^*=  H_\8^*:=- {\rm dist}(0, \cC) \in (-\8,0]\;.   
 \end{equation}
Here, $0$ is the origin in $\Z^+\times \Z^d$, and 
dist is the chemical ``distance'' given, for $s\leq t$, by 
dist$((s,x),(t,y))= \inf \{\sum_{s < u \leq t} \h(u, x_u)\}$ 
where the infimum is taken over oriented nearest neighbor 
paths $((u, x_u);s < u \leq t)$ with $x_s=x, x_t=y$. 
We note that the convergence 
$$ \mu_n^\b[H_n-n\lm'(\b)]  \longrightarrow_n  (\ln W_\8)'(\b)$$
in (\ref{eq:flenergy})
parallels that of (\ref{eq:888}), in the sense that $ \mu_n^\b(H_n)$
and $H_n^*$, which relates via  (\ref{eq:887}), both have
order one fluctuations.

As a related remark, let us recall the local limit theorem of Sinai 
\cite{sinai}. Deep inside the region $U_1$, 
$$
P\left[ \exp\{\b H_n( \om)\} \Big\vert \om_n=x\right]
= W_\8 \times \ W_\8   \circ  \theta_{n,x}^\leftarrow + R_{n,x}
$$ 
where $\theta_{n,x}^\leftarrow$ is given by
$\theta_{n,x}^\leftarrow (\h(\cdot, \cdot)) : (u,y) \mapsto
\h(n-u, x+y)$, and the error term $ R_{n,x} \to 0$ in $L^1$ uniformly
in $x: |x|\leq An^{1/2}$. The local limit theorem 
parallels
%%is the counterpart
the following observation in the percolation model
$$
H_n^{*,x} \st{\rm def}{=} \inf \{ H_n (\om) ; \om_0=0, \om_n=x\}
= - {\rm dist}(0, \cC) - {\rm dist}((n,x), \cC) +o_Q(1)
$$
for $x$ not too large.

%%%\end{example}
\medskip

{\bf Acknowledgements}: We thank Geoffrey R. Grimmett and Herbert Spohn
for stimulating conversations
and for indicating the reference \cite{gh}.

\small
%%%%%%%%%%%%%%%%%%

\end{document}